\documentclass[12pt]{article}
\usepackage{graphicx}
\usepackage{amsmath}
\usepackage{amsthm}
\usepackage{graphics}
\usepackage{epsfig}
\usepackage{amssymb}

\usepackage{xcolor}
 \usepackage{setspace}
\usepackage{extarrows}

\usepackage[authoryear]{natbib}
\usepackage{subfig}

\usepackage{pdflscape}
\usepackage{subfig}
\usepackage{enumitem}
\allowdisplaybreaks

\newcommand{\be}{\begin{eqnarray}}
\newcommand{\ee}{\end{eqnarray}}
\newcommand{\bea}{\begin{eqnarray*}}
\newcommand{\eea}{\end{eqnarray*}}

\newcommand{\lb}{\left(}
\newcommand{\rb}{\right)}

\newcommand{\N}{\mathbb{N}}
\newcommand{\Var}{\operatorname{Var}}
\newcommand{\R}{\mathbb{R}}

\newcommand{\E}{\mathbb{E}}

\newtheorem{theorem}{Theorem}[section]

\newtheorem{remark}{Remark}[section]

\begin{document}
\title{Likelihood ratio tests for many groups  in high dimensions}
\author{
{\small Holger Dette, Nina D\"ornemann} \\
{\small Fakult\"at f\"ur Mathematik} \\
{\small Ruhr-Universit\"at Bochum} \\
{\small 44799 Bochum, Germany} \\
}
\maketitle
\date{}
\begin{abstract}
In this paper we investigate the asymptotic distribution of likelihood ratio tests in models with several groups, when the number of groups converges with the dimension and sample size to infinity.
We derive   central limit theorems for the logarithm of various test statistics and compare our results with the approximations obtained from a central limit theorem where  
 the number of groups  is  fixed.

\end{abstract}
Keywords:  likelihood ratio test, high-dimensional inference

AMS subject classification:  62H15, 62H10

\section{Introduction}
\label{sec1}
\def\theequation{1.\arabic{equation}}
\setcounter{equation}{0}

Classical  multivariate analysis tools as can be found in the text books of \cite{muirhead1982} or \cite{anderson2003}
are developed under the paradigm  that the dimension is substantially smaller than the sample size  and do not
yield to a reliable statistical inference if this assumption is not satisfied.
Because modern datasets, as they  occur in biostatistics, wireless communications and finance,
are high-dimensional
[see, e.g., \cite{Fan2006}, \cite{Johnstone2006} and references therein]   there exists an enormous amount of literature
developing statistical methods in the case where the dimension of the data is of comparable size (or even larger) than the sample size.
Many authors have worked on this problem and a large part of the literature investigates the asymptotic properties of ``classical''  test procedures under the assumption that the dimension $p$ is proportional to the sample size $n$
[see \cite{ledoitwolf2002,fujhimwak2004,birkedette2005,schott2007,baietal2009,chenzhangzhong2010,kakiiwas2008,jiangyang2013,jiangetal2013,liqin2014,wang2014,jiangqi2015,hyodoetal2015,bao2017,yamadaetal2017,yanpan2017,hancheliu2017,chejia2018,chenliu2018,boddetpar2019} among many others].

\begin{figure}[t]
    \centering
    \caption{\it  Simulated $p$-values of the likelihood ratio test for the hypothesis of a block diagonal structure  of $40$ blocks
    of equal size $18$  in a $p=720$-dimensional
 normal distributed vector (sample size  $800$). Left panel: asymptotic level $\alpha$ test
      considering the number of groups as fixed; right panel: asymptotic level $\alpha$ test derived in this paper. }
      {\includegraphics[width=0.25\columnwidth]{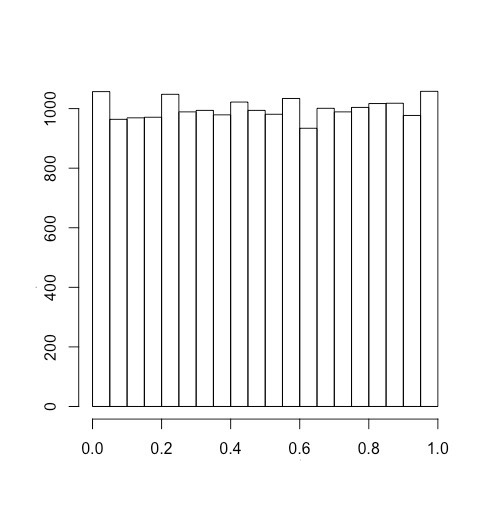}}
        \qquad
 {\includegraphics[width=0.25\columnwidth]{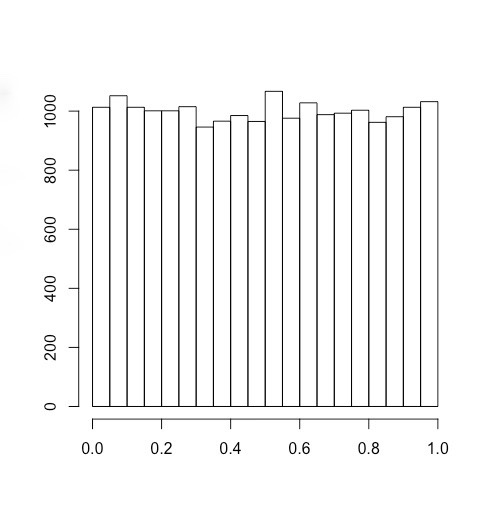}}
    \label{fig11}
\end{figure}

In the  case  $p < n$ likelihood ratio tests are still well defined and it is shown in many papers that the
asymptotic theory under the assumption  $\lim_{n,p \to \infty}  p/n =y \in (0,1)$  yields a substantially better approximation for
the nominal level of corresponding tests  as classical asymptotic considerations keeping the dimension $p$ fixed.
In this paper we continue  this discussion and investigate approximations  for the likelihood ratio test statistics
in cases  where  high-dimensional inference has to be performed for a large number of groups.
 Consider, for example, the problem of testing if the covariance matrix of  a $p$-dimensional
 normal distributed vector has  a block diagonal structure with $q$ blocks.
 In  Figure \ref{fig11} we show the $p$-values of the corresponding  likelihood ratio test
 [see \cite{wilks1935}]  under the null hypothesis with $40$ blocks of  size $18$ (thus the total dimension is $720$)
 and sample size   $800$. The components of all vectors are independent identically standard  normal distributed and thus the null hypothesis 
 is obviously satisfied.
 The left panel shows the simulated  $p$-values  (based on $20000$ simulation runs) using the
 approximation provided by  \cite{jiangyang2013} considering the number of blocks  as fixed, while the right panel
 shows the $p$-values using an approximation with $q \to \infty$ as derived in this paper.  The
 two figures look very similarly.

\begin{figure}[t]
    \centering
    \caption{\it  Simulated $p$-values of the likelihood ratio test for the equality of $q$ $p$-dimensional normal distributions, where
     $n_i=80,~p=40,~q= 300$. Left:   asymptotic level $\alpha$ test  ($p,n_{i}\to \infty$)
      considering the number of groups as fixed; right panel: asymptotic level $\alpha$ test derived in this paper for $p,n_{i},q \to \infty$ }
      {\includegraphics[width=0.25\columnwidth]{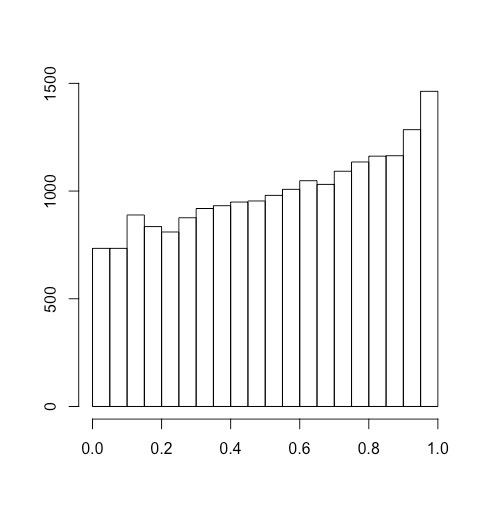}}
        \qquad
 {\includegraphics[width=0.25\columnwidth]{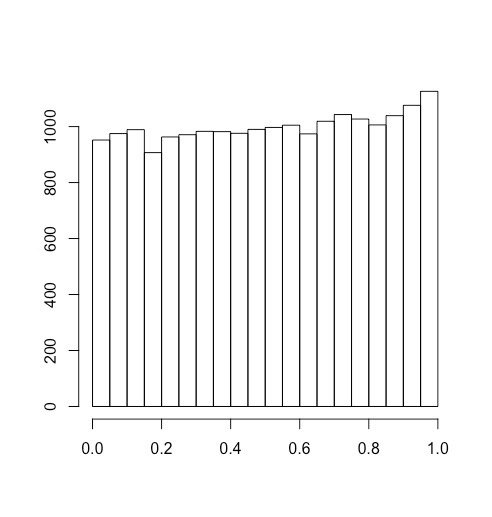}}
    \label{fig12}
\end{figure}

 On the other hand, in Figure \ref{fig12}  we show the $p$-values of the likelihood ratio test for
testing the  equality of  $q$ normal distributions  ${\cal N} (\mu_{1},\Sigma_{1}), \ldots , {\cal N} (\mu_{q},\Sigma_{q})$ and the null hypothesis $\mu_i = \mu_j; \ \Sigma_i = \Sigma_j \ (i,j=1,\ldots,q)$
[see  \cite{wilks1946}].
The sample size in each group  is $n_{i} = 50$, ($i=1,\ldots , q$), the dimension is $p=40$ and $q=300$ different groups
are considered (again all components of all vectors are independent identically standard normal distributed and $20000$ simulation runs have been performed).
 The left panel of the figure  shows the results obtained using the quantiles for the  asymptotic
distribution obtained for fixed $q$ [see Theorem  3 in  \cite{jiangyang2013}] while the right one corresponds to
an asymptotic distribution derived in this paper under the assumption that $p,q,n_i \to \infty$
(see Theorem \ref{thm3} for more details). In this case we
 observe that the latter approach provides  a better  approximation of the nominal level.

The present paper is devoted to give some (partial) explanation of  observations of this type. We consider classical
testing problems in high-dimensional statistical inference, where  data can be  decomposed in $q$ groups, and investigate
the asymptotic properties  of  likelihood ratio tests for  various  hypotheses if the dimension $p$ and the number of groups $q$  converge to infinity with
increasing sample size. In all cases  we establish the  asymptotic normality of the log-likelihood ratio after appropriate standardization. 

The work, which is most similar in spirit  to our paper is  the paper of \cite{jiangyang2013}, who considered the corresponding problems for a fixed
number of  groups.  In contrast to these authors, who used the fact that the moment generating function
of the log-likelihood ratio statistic  can essentially be expressed as a product of ratios of Gamma functions,
we  use a  central limit theorem for sums of  a triangular array of independent random variables (see Theorem \ref{clt} in Section \ref{sec5})
to establish asymptotic normality.
This approach is also applicable for other high-dimensional problems. As an example, we revisit
the problem of testing a linear hypothesis about regression coefficients as considered in
 \cite{baietal2013}. These authors showed the asymptotic normality of the (standardized) log-likelihood ratio test statistic by
using recent  results about  linear spectral statistics of large dimensional $F$-matrices.  With our
approach we are able to extend
their result and also provide a more handy  representation of the asymptotic bias.

\section{One sample problems}
\label{sec2}
\def\theequation{2.\arabic{equation}}
\setcounter{equation}{0}

\subsection{Testing for independence}  \label{sec21}

A very prominent problem in high-dimensional data analysis is the problem of testing for the independence of  sub-vectors  of a  multivariate
normal distribution. To be precise, let  $X\sim \mathcal{N}(\mu,\Sigma)$ denote a  $p$-dimensional normal distributed vector with mean
 $\mu\in\mathbb{R}^{p}$  and positive definite variance  $\Sigma\in\mathbb{R}^{p \times p}$  and assume that  $X$ is decomposed
 as
	\begin{align*}
		X = \big (  {X^{(1)}}^{\top}	, 	\ldots ,
		{X^{(q_n)} }^{\top}
		\big ) ^{\top}
		\end{align*}	
		where $X^{(i)}$ are vectors of dimension  $p_i$  ($i=1,\ldots ,q_n$)  such that $\sum_{i=1}^{q_{n}} p_{i}=p$. Let
		\begin{align} \label{22}
		\Sigma = \begin{pmatrix}
			\Sigma_{11} & \Sigma_{12} & \ldots  & \Sigma_{1q_n} \\
			\Sigma_{21} & \Sigma_{22} & \ldots  & \Sigma_{2q_n} \\
			\vdots & \vdots & ~ & \vdots\\
			\Sigma_{q_n 1} & \Sigma_{q_n 2} & \ldots  & \Sigma_{q_n q_n}
		\end{pmatrix}.
	\end{align}	
		denote the corresponding decomposition of the covariance matrix, where $\Sigma_{ij} := \operatorname{Cov}(X^{(i)},X^{(j)})$. The hypotheses of independent sub-vectors is formulated as
		\begin{align} \label{21}
		H_0: ~ \Sigma_{ij}  =0 \textnormal{~~for all~} i\neq j.
	\end{align}	
	Several authors have developed tests for the hypothesis \eqref{21} [see \cite{jiangetal2013,hyodoetal2015,bao2017,jiangqi2015,yamadaetal2017,chenliu2018,boddetpar2019} among others], and in this section we focus on the
	likelihood ratio test based on a sample of independent identically distributed observations
 $X_1,\ldots  ,X_n {\sim} {\cal N} (\mu, \Sigma ) $. \cite{wilks1935} showed  that the likelihood ratio statistic
 for the hypotheses \eqref{21} is given by
	\begin{align*}
		\Lambda_n = \frac{|\hat \Sigma |^{n/2}}{\prod\limits_{i=1}^{q_n} |\hat \Sigma_{ii}|^{n/2}},
	\end{align*}
where
	\begin{align*}
		\hat \Sigma =\frac{1}{n}\sum\limits_{k=1}^n (X_k -\overline{X})
		(X_k-\overline{X})^\top ~,
	\end{align*}
is the common estimator of the covariance matrix,	$\overline{X} ={1\over n} \sum_{k=1}^{n}X_{k}$ the sample mean
and $ \hat \Sigma_{{ij}}$ denotes the block in the $i$th row and $j$th column of  the estimate $\hat \Sigma$ corresponding to
the decomposition \eqref{22}.
The following result specifies the asymptotic distribution of the likelihood ratio test under the null hypothesis of
independent blocks, if the number of blocks  $q_{n}$ is increasing with the sample size.  A proof can be found in Section \ref{proof_thm1}.
Here and throughout this paper the  symbol $\xrightarrow{{\cal D}}$  denotes weak convergence.

\begin{theorem} \label{thm1} If  $q_{n} \to \infty $,
$p_i/n \to \lambda_i \in (0, 1)$, ${p}/{n} \to c < 1$ and $\sum _{i = 1}^\infty \lambda_i = c$,
then under the null hypothesis \eqref{21}
    \begin{align*}
     {2 \over n}   \log  \Lambda_n -  s_n \xrightarrow{{\cal D}} \mathcal{N}(0, \sigma^2)
    \end{align*}
   where
    \begin{align}
    \label{b1}
        \sigma^2 &= 2 \log \Big ((1-c)^{-1} \prod
        \limits_{i = 1}^\infty (1 - \lambda_i)\Big ) \\
        \label{b2}
        s_n &= \sum \limits_{i = 1}^{q_n} (n - p_i - 1) \log \big( 1 - \tfrac{p_i}{n - 1} \big )   - (n - p - 1)\log \big (1 - \tfrac{p}{n - 1}  \big ) - \frac{\sigma^2}{4}.
    \end{align}

\end{theorem}

    \begin{remark}  \label{rem1} {\rm
    ~~
   \begin{itemize}
    \item[(a)]   Theorem \ref{thm1}  provides an asymptotic level $\alpha$ test for the null hypothesis \eqref{21} by rejecting $H_{0}$, whenever
       \begin{align}
    \label{t1} 
			- \frac{2}{n} \log \Lambda_n > \sigma_{n,q_n} u_{1-\alpha} - s_n~,
		\end{align}			        
        where $u_{1-\alpha}$ denotes the $(1-\alpha)$-quantile of the standard normal distribution and 
        \begin{eqnarray} \label{sq}
\sigma_{n,q}^{2}  &=& 2 \log \Big ( \frac {\prod^q_{i=1}(1- \frac {p_i}{n })}{(1- \frac {p}{n })} \Big).
\end{eqnarray}
    \item[(b)]
\cite{jiangyang2013} derived the asymptotic distribution of the statistic $\frac {2}{n} \log \Lambda_n$ in the case, where the number of groups is fixed, that is $q_n=q$
and the dimension is proportional to the sample size. In particular they showed that under the null hypothesis
\begin{equation}\label{yy}
{  \frac {2}{n} \log \Lambda_n - s_{n,q} \over  \sigma_{n-1,q}} 
  \stackrel{\mathcal{D}}{\longrightarrow}  \mathcal{N}(0,1),
\end{equation}
where $ \sigma^2_{n-1,q}$  is defined in \eqref{sq} and
\begin{eqnarray*}
s_{n,q} &=& \sum^q_{i=1} (n- p_i - \frac {3}{2}) \log (1 - \frac {p_i}{n - 1}) - (n-p-\frac {3}{2}) \log (1 - \frac {p}{n - 1})
\end{eqnarray*}
(note that these authors use a slightly different notation). The corresponding  asymptotic level $\alpha$ test for the null hypothesis \eqref{21} rejects 
$H_{0}$, whenever
       \begin{align}
    \label{t2} 
        	- \frac{2}{n} \log \Lambda_n > \sigma_{n-1,q} u_{1-\alpha} - s_{q,n}~.
        \end{align}
It is easy to see that under the assumptions of Theorem \ref{thm1}
$$
\lim_{n \to \infty} \sigma^2_{n-1,q_{n}} = \sigma^2~,
$$
where $\sigma^2$ is defined in \eqref{b2}. Moreover, recalling the definition of $s_n$ in \eqref{b2} we obtain by a straightforward calculation
$$
\lim_{n \to \infty}  \big (
s_n - s_{n,q_n}  \big ) = \lim_{n \to \infty}  \Big (
 \frac {1}{2} \sum^{q_n}_{i=1} \log (1 - \frac {p_i}{n - 1}) - \frac {1}{2} \log (1- \frac {p}{n - 1}) - \frac {\sigma^2}{4} \Big ) =0 .
$$
These results explain  why  in   Figure \ref{fig11} the simulated $p$-values of the likelihood ratio test \eqref{t2} 
 obtained  by a central limit theorem
with $p,n \to \infty$, $q$ fixed  and the likelihood ratio test \eqref{t1}  obtained  by a central limit theorem using
 $p,n, q\to \infty$ are very similar.
\end{itemize}
}
    \end{remark}

  \subsection{Testing a linear hypothesis about regression coefficients}
  \label{sec22}
  A further problem appears  if the $p$-dimensional  (independent)  random variables  $X_1, \ldots , X_n$
  depend  linearly  on  $q$-dimensional regressors, say $z_{1},\ldots , z_{n}$.   To be precisely, assume
  $X_k \sim \mathcal{N}(\beta z_k, \Sigma)$ ($1\leq k \leq n$), where the covariance matrix $\Sigma\in\R^{p\times p}$
  is positive definite and $z_1,\ldots , z_n \in\R^q$ are known design vectors such that the matrix $(z_1,\ldots , z_n )$  has rank $q$.
Consider the decomposition
	\begin{align*}
		\beta=(\beta_1,\beta_2),
	\end{align*}
	with $p\times q_1$  and  $p\times q_2$ matrices  $\beta_1$  and  $\beta_2$, respectively,   such  that $q=q_1 + q_2$.
We are interested in the hypothesis that  the matrix  	$\beta_{1}$ coincides with  a given matrix $\beta_{01}\in\R^{p\times q_1}$, that is
	\begin{align} \label{c3}
		H_0: \beta_1=\beta_{01}.
	\end{align}
The likelihood ratio test statistic  for this hypothesis is given by
	\begin{align*}
	\Lambda_n :=  \frac{| \hat{\Sigma}| ^\frac{n}{2}}{|\hat{\Sigma}_0|^\frac{n}{2} }  ,
	\end{align*}
	where
	the  $ p \times p  $ matrices $\hat{\Sigma}$ and $\hat{\Sigma_{0}}$ are defined by
	\begin{align*}
			\hat{\Sigma} &= \frac{1}{n} \sum\limits_{i=1}^n  (  X_i - \hat{B} z_i ) ( X_i - \hat{B} z_i )  ^\top  , \\
				\hat{\Sigma}_0 &= \frac{1}{n} \sum\limits_{i=1}^n (X_i - \beta_{01} z_{i,1}  - \hat{B}_{20}z_{i,2}) (y_i - \hat{B}_{20} z_{i,2} )^\top ,
	\end{align*}
	respectively, and
	\begin{align*}
		\hat{B}  & =  \Big ( \sum\limits_{i=1}^n X_i z_i^\top  \Big )  \Big ( \sum\limits_{i=1}^n z_i z_i^\top \Big )^{-1}, \\
	\hat{B}_{20} & =  \Big ( \sum\limits_{i=1}^n (X_i - \beta_{01} z_{i,1})  z_{i,2}^\top \Big ) \Big ( \sum\limits_{i=1}^n z_{i,2} z_{i,2}^\top  \Big )^{-1}
	\end{align*}
	are
	 the maximum likelihood estimators of $\beta$ under the null hypothesis and alternative, respectively
	[see  \cite{sugiura1969} or   \cite{anderson2003}].
Here  we use the  partition  of the vector  $z_i^\top = (z_{i,1}^\top, z_{i,2}^\top)$ in vectors  $z_{i,1}^\top$ and
	$ z_{i,2}^\top$ of dimension $q_1$ and $q_2$, respectively.
  In the following theorem, we present the asymptotic null distribution of the likelihood ratio test statistic for  a general linear hypothesis
  \eqref{c3}  in a high-dimensional regression model, where the dimensions $p=p_{n}$, $q=q_{n}$, $q_{1} = q_{1,n}$ and $q_{2} = q_{2,n}$
  increase with the sample size.
  A part of this result, namely the  case $\frac{p_{n}}{q_{1,n}}\to y_1 \in (0,1)$,
    has been established by  \cite{baietal2013} using random matrix theory.
  In contrast to these authors  we are also able to deal with  the case $ y_1 \in (0,\infty)$. 

	\begin{theorem}\label{thm2}
	If $p\to \infty, ~ q_{1,n} \to \infty,$ $n - q_n\to \infty$, $\frac{p}{q_{1,n}}\to y_1 \in (0,\infty)$  and $\frac{p}{n-q_{n}} \to y_2 \in (0,1)$, then under the null hypothesis \eqref{c3}
	\begin{align*}
		\frac{2}{n} \log \Lambda_n - s_n \xrightarrow{{\cal D}} \mathcal{N}(0,\sigma^2),
	\end{align*}
	where
	\begin{align}
	\sigma^2:= &2 \Big\{ \log \Big ( \frac{1}{1-y_2} \Big )  - \log \Big (  \frac{y_1 + y_2}{y_1 + y_2 - y_1 y_2} \Big  ) \Big \}, ~\label{c1} \\
	s_n :=& ( n - q_{2,n} - 1) \log \Big ( \frac{ n - q_{n} - 1} { n - q_{2,n} - 1}  \Big )
	+ q_{1,n} \log  \Big ( \frac{n - q_{n}- p - 1} {n -  q_{n} - 1} \Big )
	\nonumber\\
	&   + ( n - q_{2,n} - p - 1) \log  \Big ( \frac{n - q_{2,n} - p - 1}{ n - q_{n} - p - 1}  \Big )  + \frac{\sigma^2}{4}. \label{c2}
	\end{align}
	\end{theorem}		
	
	\begin{remark}\label{rem22}
	{\rm
	\cite{baietal2013} considered the  testing problem \eqref{c3}
	in a similar high-dimensional framework.
		Note that the authors use the negative log likelihood ratio test statistic $-\log \Lambda_{n}$, while Theorem \ref{thm2} is formulated for $\log \Lambda_{n}$.
		They  made use of  recent  results about
	linear spectral statistics of large dimensional $F$-matrices
	and require a more restrictive condition on the ratio $ \frac {p_{n}}{q_{1,n} }$ to apply this theory, that is
	 $\lim_{n\to \infty } \frac {p_{n}}{q_{1,n} } = y_1\in (0,1)$.
	To be more  precise, \cite{baietal2013}   proved that under the null hypothesis \eqref{c3}
	\begin{align} \label{baiclt}
	v(f)^{-\frac{1}{2}} \big (  - \frac{2}{n}\log \Lambda_n - p F_{y_{n_1},y_{n_2}}(f) - m(f) \big )  \stackrel{\mathcal{D}}{\to} \mathcal{N}(0,1),
	\end{align}
	whenever $p\to \infty, ~ q_{1,n} \to \infty, ~n - q_n\to \infty$, $\frac{p}{q_{1,n}}\to y_1 \in (0,1)$  and $\frac{p}{n-q_{n}} \to y_2 \in (0,1)$.
	For an explicit definition of the expression $m(f),~v(f)$ and $F_{y_{n_1},y_{n_2}}$, we refer the reader to formulas
	(26), (27) and (29) in their paper.
	Theorem \ref{thm2}  extends the result of \cite{baietal2013}  to the case where $y_{1} \geq 1$ and provides a simpler representation of the bias. Moreover,
	we have checked numerically that the standardizing terms in  the central limit theorem stated in \eqref{baiclt}  and
	Theorem \ref{thm2} behave similarly. \\
Consequently,  the likelihood ratio tests based on the asymptotic distribution of Theorem \ref{thm2} and Theorem  3.1 of  \cite{baietal2013} 
	have very similar properties. Numerical results, which confirm this observation are not displayed for the the sake of brevity.}
	\end{remark}
 	
	  \section{Some $q$-sample problems}
	  \label{sec3}
\def\theequation{3.\arabic{equation}}
\setcounter{equation}{0}
	
 In this section we consider the comparison of $q$  normal distributions 
	  $\mathcal{N}(\mu_1, \Sigma_1), \ldots ,  \mathcal{N}(\mu_q, \Sigma_q)$ with mean vectors
   $\mu_{1}, \ldots , \mu_{q}\in\R^p$ and covariance matrices
   $\Sigma_{1} , \ldots , \Sigma_{q}
   \in\R^{p\times p}$, where for each  group
   a sample of size $n_j$ is available ($j=1,\ldots ,q$)  and the dimension and number of groups are increasing with the sample size.	
	  
  \subsection{Testing equality of several covariance matrices}
  \label{sec23}
An important assumption for  multivariate analysis of variance (MANOVA)  is that  of  equal
   covariances in the different groups. Thus we are interested in a  test of the hypothesis
 	\begin{align} \label{c4}
 		H_0: \Sigma_1 = \ldots  = \Sigma_q.
 	\end{align}
	This problem has been considered by several authors in the context of high-dimensional inference
	[see  \cite{obrien1992}, \cite{schott2007}, \cite{srivastava2010} or  \cite{jiangyang2013} among others].

	In this section  we add to this line of literature and
	investigate the asymptotic distribution of  the likelihood ratio test based on samples of independent distributed observations
	$X_{ji} \sim \mathcal{N}(\mu_j, \Sigma_j),$ $1\leq i \leq n_j$, $1\leq j \leq q$, when the number of groups is large, i.e. $q \to \infty$. To be precise,   let
$n = \sum_{j=1}^q n_j $ 		be the total sample size,
then the test statistic of the likelihood ratio test  for the hypothesis \eqref{c4}  was derived by  \cite{wilks1932} and is  given by
	\begin{align} \label{h1}
		\Lambda_{n,1} = \frac{\prod\limits_{j=1}^q |A_j / n_j |^{\frac{1}{2}n_j }}{|A/n|^{\frac{1}{2}n }},
	\end{align}
	where  the $p \times p$  matrices   $A_{j}$ and $A$ are defined as
	\begin{align}
		A_j &= \sum_{k=1}^{n_j} (X_{jk} - \bar{X}_j)(X_{jk} - \bar{X}_j)^\top ~,~~
		A = \sum_{j=1}^q A_j.\label{c5}
	\end{align}
	As proposed by \cite{bartlett1937} we consider the modified likelihood ratio test statistic
	\begin{align*}
	\tilde{\Lambda}_{n,1} =  \frac{\prod\limits_{j=1}^q |A_j/(n_j-1)|^{\frac{1}{2} (n_j - 1)}}{|A/(n-q)|^{\frac{1}{2}(n - q )  }},
	\end{align*}
where each sample size $n_j$ is substituted by its degree of freedom. Our next result deals with asymptotic distribution of the test statistic $\log \tilde{\Lambda}_{n,1}$ 
 for an increasing dimension and an increasing number of groups.

  \begin{theorem} \label{thm4}
  Let $n_j+1>p$ for all $1\leq j \leq q$,  $p\to \infty$, $q=o(\sqrt{p (n-q)})$, $p=o(n-q)$, assume that
$$
\limsup_{n\to\infty} \sup\limits_{1\leq j \leq q} \frac{p}{n_j-1} <1  ~.
$$
and that  
$$
  \sigma_n^2   =  \frac{1}{2} \sum\limits_{j=1}^q \frac{(n_j - 1)^2}{p(n -q )} \log \Big ( \frac{n_j - 1 }{n_j - p - 1} \Big )  - \frac{1}{2},
$$
converges with a positive limit, say $\sigma^{2} >0$. 
 Then, under the null hypothesis \eqref{c4},
  \begin{align*}
  \frac{\log \tilde{\Lambda}_{n,1} - \tilde{s}_n }{\sqrt{p(n-q)}} \to \mathcal{N}(0,\sigma^2),
  \end{align*}
  where
  \begin{align*}
  \tilde{s}_n
    &=  \sum\limits_{j=1}^q \frac{n_j - 1}{2} \Big\{ \Big ( n_j - \frac{3}{2} \Big )   \log  \Big (\frac{n_j - 2}{n_j - p - 2} \Big )   - p \log \lb \frac{n_j - 1}{n_j - p - 2} \rb \Big\} \\
    & ~~~ - \frac{n - q}{2}  (n - q - p) \log  \Big ( \frac{n - q}{n - q - p}  \Big ),    
  \end{align*}

  \end{theorem}

 \begin{remark} \label{rem4}  ~~~\\
 {\rm
 \begin{itemize}
\item[(a)] Note that under the assumptions  of Theorem \ref{thm4} the asymptotic distributions of $ \log\Lambda_{n,1} $  and $ \log \tilde{\Lambda}_{n,1}$  are not
 identical. In fact we have
\begin{align*}
	\log\Lambda_{n,1} - \log \tilde{\Lambda}_{n,1}
	&= \frac{1}{2} \sum_{j=1}^q p n_j \log \Big ( \frac{n_j - 1}{n_j} \Big ) - \frac{1}{2}pn \log \Big ( \frac{n-q}{n} \Big )
	- \sum\limits_{j=1}^q \frac{1}{2} p \log \Big ( \frac{n_j - 1}{n -q} \Big )   ,
\end{align*}	
and in general this is not of order $o(\sqrt{p(n-q)})$  (consider for example the case $n_{j}=2p+1$, $q=o(p^{2})$, $q \to \infty$).
\item[(b)]\cite{jiangyang2013} determined the asymptotic distribution of the statistic $\log \tilde{\Lambda}_{n,1}$ for a fixed number $q$ assuming that the limit $p/n$ does not vanish. In particular, they showed that under the null hypothesis \eqref{c4}
	\begin{align} \label{c8}
		\frac{\log \tilde{\Lambda}_{n,1} - \mu_n }{(n-q)\sigma^{(1)}_n} \stackrel{\mathcal{D}}{\to} \mathcal{N}(0,1)
	\end{align}
	as $n_j,n,p\to\infty$, where the asymptotic bias and variance are given by
	\begin{align*}
		\mu_n &= \frac{1}{4} \Big \{ (n - q) (2n - 2p - 2q - 1) \log \Big( 1 - \frac{p}{n-q} \Big) \\
		&~~~
		- \sum\limits_{j=1}^q (n_j - 1) (2n_j - 2p -3) \log \Big( 1 - \frac{p}{n_j - 1} \Big) \Big\}, \\
		(\sigma_n^{(1)})^2 & = \frac{1}{2} \Big\{ \log \Big( 1 - \frac{p}{n-q} \Big) - \sum\limits_{j=1}^q \Big( \frac{n_j - 1}{n - q} \Big)^2 \log \Big( 1 - \frac{p}{n_j - 1} \Big) \Big\},
	\end{align*}
respectively.	As $q$ being fixed, the authors assumed for their result that $p/n_j \to y_j \in (0,1]$ for all $j=1,\ldots ,q$ and $\min_j n_j > p + 1$. Note that the order  of  standardization in   Theorem \ref{thm4} is different than in \eqref{c8}. The standardization is of order $\sqrt{p(n-q)}$ which is, under the assumptions of Theorem \ref{thm4}, substantially smaller than $n-q$ as used in \eqref{c8}. 
Comparing  the variance $\sigma^2$ in Theorem \ref{thm4} with an adjusted version of $(\sigma_n^{(1)})^2$ (such that the different standardizations are canceled out)
yields under the assumptions of Theorem \ref{thm4}
	\begin{align*}
	\sigma^2 - \frac{n-q}{p}(\sigma_n^{(1)})^2 = - \frac{1}{2} - \frac{1}{2} \frac{n - q}{p} \log \Big(  1 - \frac{p}{n-q} \Big) = o(1).
	\end{align*}
On the other hand   the difference of the means is given by 
	\begin{align*}
	\tilde{s}_n - \mu_n =& \frac{n-q}{4} \log \Big (  \frac{ n - q - p}{n - q}  \Big ) + \sum\limits_{j=1}^q \frac{n_j - 1}{2} \Big\{ p \log \Big(  \frac{n_j - p - 2}{n_j - p - 1}  \Big )  \\
	& +  \Big ( n_j - \frac{3}{2}  \Big)  \Big [ \log  \Big (  \frac{n_j - 2}{n_j - p -2}  \Big )  - \log  \Big (  \frac{n_j - 1}{n_j - p - 1} \Big)  \Big ] \Big\}
	\end{align*}
	Note that the first summand divided by the standardization $\sqrt{p(n-q)}$ vanishes under the assumptions of Theorem \ref{thm4}, while the other 
	terms  give a notable contribution to the expected value.  Thus we expect that the corresponding likelihood ratio tests behave differently, if the number of groups is large.
	\end{itemize}
	}
	\end{remark}	
	
  \subsection{Testing equality of several normal distributions}
  \label{sec24}
 	We consider the same setting as in Section \ref{sec23} but this time want to test whether $q$ normal distributions are identical, that is,
 	\begin{align} \label{c6}
 		H_0: \mu_1=\ldots =\mu_q, ~\Sigma_1 = \ldots   = \Sigma_q.
 	\end{align}
	The  test statistic of the  likelihood ratio  test for the hypothesis \eqref{c6} is given by
 	\begin{align} \label{c7}
 	\Lambda_n = \frac{\prod\limits_{j=1}^q |A_j / n_j| ^ {\frac{1}{2}n_j }}{|B/n|^{\frac{1}{2}n}},
 	\end{align}
	where the $p\times p$ matrix $A$   is defined in \eqref{c5} and
 	\begin{align*}
 	B = \sum\limits_{j=1}^q \sum\limits_{i=1}^{n_j} (X_{ji} - \bar{X})(X_{ji} - \bar{X})^\top
 	= A+ \sum\limits_{j=1}^q n_j (\bar{X}_j - \bar{X}) (\bar{X}_j - \bar{X})^\top.
 	\end{align*}
 Note that $\Lambda_{n}$ is the product of the likelihood ratio statistic  $\Lambda_{n,1}$ in \eqref{h1} for testing equality of covariance matrices
 and the likelihood ratio test statistic for testing equality of the means [see \cite{yaobaizheng2015}, \cite{gregory2015}].
 Several authors dealt with testing $H_0$ [e.g. see \cite{wilks1946}, \cite{jiangyang2013}].
	The following result specifies the asymptotic distribution of the  statistic $ \log \Lambda_n$
 for increasing dimension and an increasing number of groups.

	\begin{theorem} \label{thm3}
	Let  {$n_j>p+1$ for all $1\leq j \leq q$}, {$p\to \infty$}, $q=o(\sqrt{p n})$, $p=o(n)$, assume that 
	$$\limsup\limits_{n\to\infty} \sup\limits_{1\leq j \leq q} \frac{p}{n_j} < 1  $$
and that 
$$
 \tilde \sigma_n^2  = \frac{1}{2}\sum\limits_{j=1}^{q} \frac{n_j^2}{p n}
		\log\Big ( \frac{n_j-1}{n_j-p-1}	\Big ) - \frac{1}{2}
		$$
		converges to a positive limit, say  $\sigma^{2 } > 0.$
Then, under the null hypothesis 
	 \eqref{c6} we have
	 	\begin{align*}
	\frac{\log \Lambda_n - \tilde s_n }{\sqrt{pn}} \stackrel{\mathcal{D}}{\to} \mathcal{N}(0,\sigma^2). 
	\end{align*}
	where
		\begin{align*}
	 \tilde s_n :=&   		\sum\limits_{j=1}^q  {\frac{n_j}{2}}\Big \{ \Big (  n_j - p - \frac{3}{2}  \Big )  \log  \Big (   \frac{n_j - 2}{n_j - p - 2}  \Big )  + p \log  \Big ( \frac{n_j - 2}{n - q}  \Big )  - p \log (n_j)  \Big \}
		 \\&+ \frac{n}{2} \Big\{ p \log  \Big (  \frac{\frac{1}{2}n - 1}{\frac{1}{2}n + \frac{1}{2}q - \frac{3}{2}}  \Big )
		  + (n - q - p) \log  \Big (  \frac{n - q - p}{n - q}   \Big ) + p \log (n)  \\
		 &+  \Big (  \frac{1}{2}n - p - 1 \Big )  \log \Big ( \frac{\frac{1}{2}n - 1}{\frac{1}{2}n - p - 1}  \Big )   +  \Big (  \frac{1}{2}n - p + \frac{1}{2}q - \frac{3}{2}  \Big )
		 \log  \Big (  \frac{\frac{1}{2}n - p + \frac{1}{2}q - \frac{3}{2}}{\frac{1}{2}n + \frac{1}{2}q - \frac{3}{2}} \Big ) \Big\} 
	\end{align*}
	\end{theorem}

  \begin{figure}[!ht]
    \centering
    \caption{
    \it  Simulated $p$-values of the likelihood ratio test for the hypothesis 
   \eqref{c6} under the null hypothesis ($n_j=200,$ $p=100$, $q=50$).  
     Left panel: asymptotic level $\alpha$ test   considering the number of groups as fixed; right panel: asymptotic level $\alpha$ test derived from  Theorem \ref{thm3}. }
    \subfloat[]{\includegraphics[width=0.3\columnwidth]{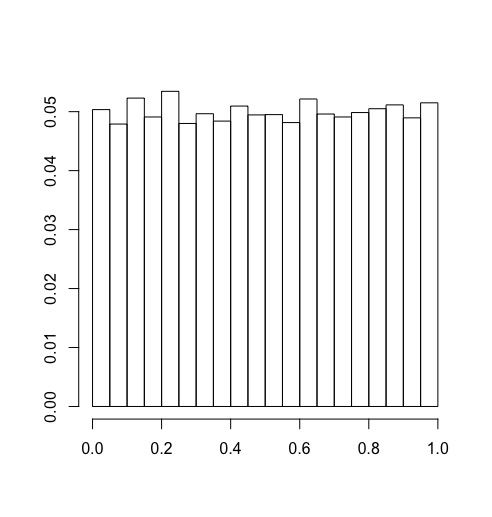}}
        \qquad
     \subfloat[]{\includegraphics[width=0.3\columnwidth]{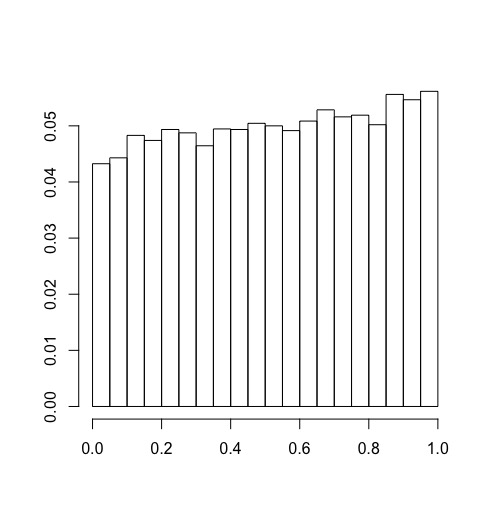}}
    \label{fig2}
\end{figure}

\begin{figure}[!ht]
    \centering
    \caption{   \it  Simulated $p$-values of the likelihood ratio test for the hypothesis 
   \eqref{c6}  under the null hypothesis, where $n_j=80,$ $p=40$ and  $q=100$  (left column),   $q=200$  (middle column), $q=300$  (right column).
Upper part: asymptotic level $\alpha$ test   considering the number of groups as fixed; lower part: asymptotic level $\alpha$ test derived from  Theorem \ref{thm3}.   }
         \includegraphics[width=0.3\columnwidth]{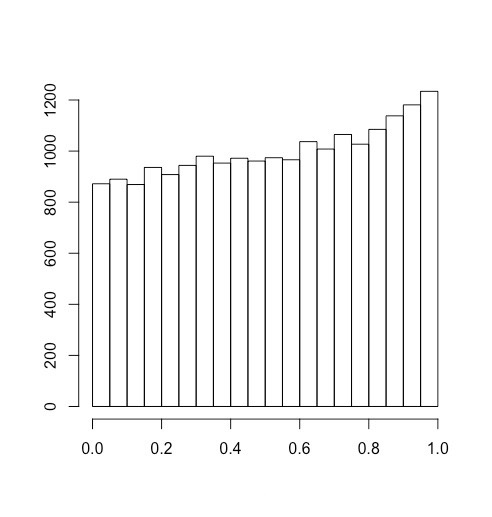}
         \includegraphics[width=0.3\columnwidth]{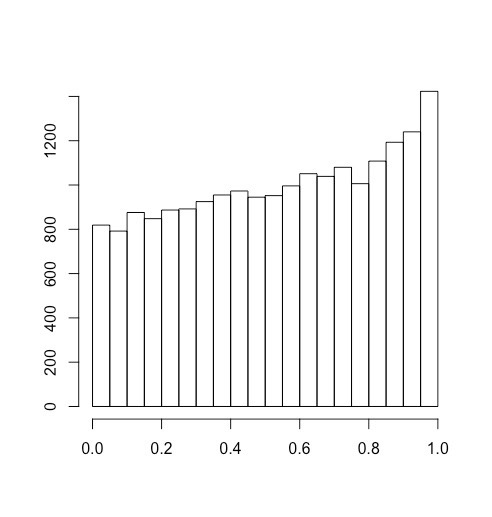}
         \includegraphics[width=0.3\columnwidth]{pval_yj_p40_q300_m80}
\includegraphics[width=0.3\columnwidth]{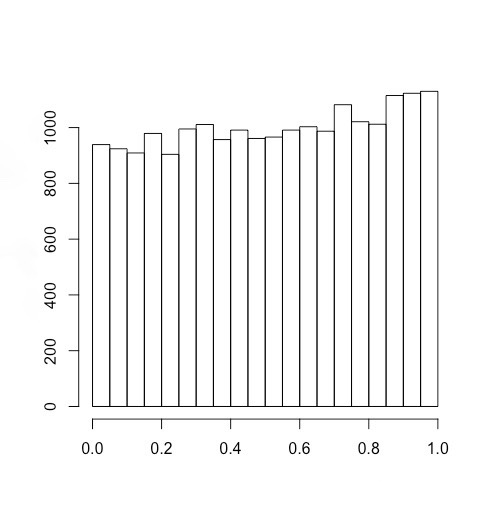}
\includegraphics[width=0.3\columnwidth]{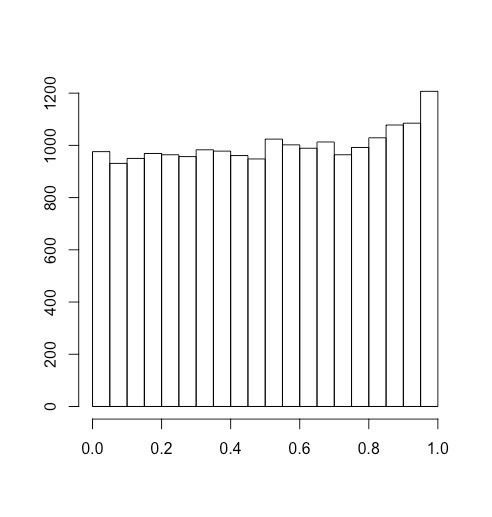}
 \includegraphics[width=0.3\columnwidth]{pval_wir_p40_q300_m80}
    \label{fig1}
\end{figure}

 \begin{remark}  {\rm 
The asymptotic distribution of the statistic $\Lambda_{n}$ in the case where $q$ is fixed was determined in Theorem 3 of
\cite{jiangyang2013} who  showed that under the null hypothesis \eqref{c6}
  	\begin{align} \label{jjthm3}
  		\frac{\log \Lambda_n -  \tilde \mu_n}{n   \tilde \sigma_n^{(1)}} \stackrel{\mathcal{D}}{\to} \mathcal{N}(0,1),
  	\end{align}
  if
$\lim\limits_{n\to\infty} \frac{p}{n_j} = y_j \in (0,1]$  ($1\leq j \leq q$). Here the asymptotic bias and
variance are given by
  	\begin{align*}
  		 \tilde \mu_n &= \frac{1}{4} \Big \{ -2qp - \sum\limits_{j=1}^q y_j - n(2p - 2n + 3)\log\Big ( 1-\frac{p}{n} \Big  ) \\
		& ~~+ \sum\limits_{j=1}^q n_j (2p - 2n_j + 3)
  		\log \Big (  1 - \frac{p}{n_j - 1}\Big )  \Big  \}, \\
  		(  \tilde \sigma_n^{(1)})^2 &= \frac{1}{2} \Big \{\log \Big ( 1 - \frac{p}{n} \Big )  - \sum\limits_{j=1}^q \frac{n_j^2}{n^2} \log \Big  ( 1 - \frac{p}{n_j - 1} \Big  )  \Big \}
  	\end{align*}
	(note that these authors use a slightly  different notation).
	 It is important to note that  the  orders in the standardizations in both results  are  different. While the standardizing factor
	 in \eqref{jjthm3} is of order $n$, it is of order $\sqrt{pn} = o(n) $ in Theorem \ref{thm3}.
	 
Similarly, as in Remark  \ref{rem4} it can be shown that 
	 $$
	 \tilde \sigma_{n}^{2} -\frac{n}p (  \tilde \sigma_n^{(1)})^2 	=o(1)
	  $$
	  under the assumptions of Theorem \ref{thm3}, while in general the difference $ \tilde s_{n}-  \tilde \mu_{n}$ is not of order $ o(\sqrt{pn} )$  (consider, for example, the case $q=p$, $n_{j} =2p+1$, $n=p(2p+1)$, 
	  $p \to \infty$). 	 
	Based on these observations we expect differences in the likelihood ratio test, 
	 if the quantiles from the normal approximation for fixed $q$ as derived  \cite{jiangyang2013} or the quantiles from Theorem \ref{thm3} are used as 
	 critical values. This is illustrated in Figure \ref{fig2} and  \ref{fig1}, where we display the simulated  $p$-values for the tests under the null hypothesis.
	 In Figure \ref{fig2} we consider the  case $n_j=200$, $p=100$ and a relatively small number of groups $q=50$. We observe that both approximations
	 yield histograms close to  the expected uniform distribution. On the other hand in  Figure \ref{fig1} we consider the cases
$n_j=80$, $p=40,$ and $q=100,200$ and $300$  and we observe larger differences in both approximations.  In particular the 
 critical values derived from Theorem \ref{thm3}  yield a likelihood ratio test for the hypothesis \eqref{c4}  with a
better performance than the test using the quantiles from  fixed $q$ asymptotics.
	}
  \end{remark}

\section{Some proofs}  \label{sec5}
\def\theequation{A.\arabic{equation}}
\def\thetheorem{A.\arabic{theorem}}
\setcounter{equation}{0}

In this section we  present proofs of our results, where we restrict ourselves to the proofs of Theorem \ref{thm1} and \ref{thm2}. 
The other statements are shown by similar arguments, which are omitted for the sake of brevity. 

\subsection{A central limit theorem}

We begin stating  a central limit theorem, which is used in the proofs of Theorem \ref{thm1} - \ref{thm3}.
We make extensive use of a central limit theorem for triangular array of independent random variables, which follows
by similar arguments as given in \cite{tomdet2019}. Therefore the proof is omitted.

\begin{theorem} \label{clt}
Let $(T_n)_{n\in\mathbb{N}}$ be a sequence of finite sets $\{ (X_n(i))_{i\in T_n} | n \in \mathbb{N} \} $
denote  an array of random variables  and $\{ (g_n(i))_{i\in T_n} | n \in \mathbb{N} \} $
an array of weights  satisfying the following conditions:
	\begin{enumerate}[label=(\textbf{A.\arabic*})]
	\item The random variables $(X_n(i))_{i\in T_n}$ are independent for all $n\in\N$. \label{A1}
	\item The random variables $X_n(i)$ are centered, that is, $\E [X_n(i)] = 0 ~ \forall i\in T_n,~n\in\mathbb{N}$. \label{A2}
	\item $ \E \left[X_n^{4} (i) \right] \leq C \E \left[ X_n^2(i)\right]^{2}$ for some universal constant $C>0$ and for all $n\in\N$. \label{A3}
	\item $\sup\limits_{i \in T_n}  g_n^2(i) \Var (X_n(i))  \to 0 \textnormal{ for } n\to \infty$. \label{A4}
	\item \label{A5} There exists a constant $\sigma^2>0$ such that
	\begin{align*}
		\sum\limits_{i\in T_n} g_n^2(i) \Var (X_n(i)) \to \sigma^2 \textnormal{ for } n \to \infty.
	\end{align*}
	\end{enumerate}
	Then  the random variable
	$$
	Z_n := \sum\limits_{i\in T_n} g_n(i) X_n(i)
	$$
	converges in distribution to a normal distribution with mean $0$ and variance $\sigma^2$.
	\end{theorem}

\subsection{Proof of Theorem \ref{thm1} } \label{proof_thm1}
  Define
	\begin{align*}
		V_n = \Lambda_{n}^{2/n} = \frac{|\hat \Sigma |}{\prod\limits_{i=1}^{q_n} |\hat \Sigma_{ii}|}
	\end{align*}
	and note that under the null hypothesis  \eqref{21} the distribution of $V_{n}$ is given by a product of independent Beta-distributions [see \cite{anderson2003}], that  is
	\begin{align*}
		V_n \stackrel{\mathcal{D}}{=} \prod\limits_{i=2}^{q_n} \prod\limits_{j=1}^{p_i} V_{ij},
	\end{align*}
	where the random variables 	
	\begin{align*}
		V_{i, j} \sim \beta((n-p_i^*-j)/2, p_i^*/2)
	\end{align*}	
 are independent and   $p_i^* = \sum 	 \limits_{l = 1}^{i - 1} p_l$. Consequently, with the notation
    \begin{align*}
      S_n:=\log V_n = \sum \limits_{i = 2}^{q_n} \sum \limits_{j = 1}^{p_i} \log(V_{i, j})
    \end{align*}
    the assertion follows from
    \begin{align} \label{b4}
        S_n - s_n \xrightarrow{{\cal D}} \mathcal{N}(0, \sigma^2)
    \end{align}
   where        $ \sigma^2 $ and $  s_n $ are defined in \eqref{b1} and \eqref{b2}, respecively.

  For a proof of this statement we use Theorem \ref{clt} and show  that  the conditions  \ref{A1}-\ref{A5} in this result are
  satisfied. We begin with a calculation of the variance of $S_{n}$ noting that the variance of logarithm of a Beta distributed
  random variable $B\sim \beta(a, b)$  is given by
  \begin{align} \label{31}
  \operatorname{Var}(\log(X)) = \psi_1(a) - \psi_1(a + b),
  \end{align}
  where  $\psi_k(x) = \frac{d^{k + 1}}{dx^{k + 1}} \log \Gamma(x)$  ($k \geq 0$) denotes the polygamma function of order $k$
  [see \cite{abrsteg1964}]. This yields
      \begin{align} \label{var}
        \operatorname{Var}(S_n) = A_{n}^{(1)}   - A_{n}^{(2)},
    \end{align}
    where
       \begin{align*}
A_{n}^{(1)}    &= \sum \limits_{i = 2}^{q_n}  \sum \limits_{j = 1}^{p_i} \psi_1((n-p_i^*-j)/2)  \\
A_{n}^{(2)}    &=   \sum \limits_{i = 2}^{q_n}  \sum \limits_{j = 1}^{p_i}  \psi_1((n  - j)/2).
      \end{align*}
Observing the expansion
  $\psi_1(z) = \frac{1}{z} + O(z^{-2})$ for the logarithmic Gamma function
  of order $1$,
  [see \cite{abrsteg1964}] we obtain for the first term
	     \begin{align}
    \nonumber
    A_{n}^{(1)}   & = \sum \limits_{i = n  - p}^{n -p_1-1} \psi_1(i/2) = \sum\limits_{i= n - p}^{n - p_1 -1} \frac{2}{i} + o(1) \\[1ex]
       & =  2 \log \Big (\frac{1 - p_1/n-1/n}{1 - p/n-1/n}\Big ) + O(n^{-1}) =  2 \log \left(\frac{1 - \lambda_1}{1 - c}\right) ( 1 +{o}(1)),
              \label{sum}
    \end{align}
    where we used the expansion
    \begin{align}\label{harmon}
		\sum\limits_{k=1}^n \frac{1}{k} = \log n + \gamma
		+ \frac{1}{2n} + O(n^{-2})
	\end{align}	
(here  $\gamma$ denotes the  Euler-Mascheroni constant).
For the second term we use the same expansion as above
and obtain
    \begin{align} \nonumber
       A_{n}^{(2)} &=      \sum \limits_{i = 2}^{q_n} \sum \limits_{j = 1}^{p_i} \psi_1((n - j)/2) = \sum \limits_{i = 2}^{q_n} \sum \limits_{j = 1}^{p_i}\Big \{ \frac{2}{n - j}  + O((n  - j)^{-2})\Big\}
       \\
 \label{a1}
      &=  \sum \limits_{i = 2}^{q_n} \sum \limits_{j = 1}^{p_i} \frac{2}{n - j}   ~+ O(n^{-1})~,
    \end{align}
    where we used for last equality that     $p_i \le n - p_1$ for sufficiently large $n$
   and the fact that  $p_{i} / n \to \lambda_{i}$ for all $i \ge 1$.
Observing the expansion for the harmonic series in \eqref{harmon}  it follows that
    \begin{align}
    A_{n}^{(2)} &=     \sum
         \limits_{i = 2}^{q_n} \sum \limits_{j = n - p_i}^{n-1} \frac{2}{j}  ~+ O(n^{-1})
         \nonumber\\[1ex]
         &= \sum\limits_{i=2}^{q_n} 2\Big\{ \log (n-1) + \frac{1}{2(n-1)}+
         O(n^{-2}) \nonumber\\[1ex]
         &{} - 2 \Big [ \log(n-p_i-1) +
         \frac{1}{2(n-p_i-1)} + O((n-p_i-1)^{-2}) \Big ] \Big\} ~+ O(n^{-1}) \nonumber\\[1ex]
        &= \sum \limits_{i = 2}^{q_n} 2 \Big\{-\log\Big (1-\frac{p_i}{n-1}
        \Big ) + \frac{p_i}{2(n-1)(n-p_i-1)}  \Big\} ~+ O(n^{-1}).  \nonumber \\
           &= \sum \limits_{i = 2}^{q_n} -2 \log\Big (1-\frac{p_i}{n-1}
        \Big ) ~+ O(n^{-1}). \nonumber \\
        & = \Big (  - 2 \sum \limits_{i = 2}^\infty \log(1 - \lambda_i)  \Big ) \big ( 1 + o(1) \big ) ~.
        \label{sum2}
    \end{align}
    Here the last equality is a consequence of the theorem of the dominated convergence [see \cite{kallenberg1997}, Theorem 1.21].
    Combining \eqref{var}, \eqref{sum} and  \eqref{sum2} finally shows
     \begin{align*}
    	\operatorname{Var}(S_n)
    	&= 2 \log \Big ( (1-c)^{-1}
    	\prod\limits_{i=1}^\infty (1-\lambda_i )\Big )
    	 + o(1)~,
	     \end{align*}
	which yields the asymptotic variance $\sigma^2$ in \eqref{b1} and assumption \ref{A5} in Theorem \ref{clt}.

Moreover, as the function $\psi_{1}$ is positive and decreasing we have
    \begin{align*}
    	\sup\limits_{i, j} \operatorname{Var}(\log(V_{i, j})) &= \sup\limits_{i,j}
    	\Big\{ \psi_1\Big(\frac{n-p_i^\star-j}{2}\Big) - \psi_1\Big(
    	\frac{n-j}{2} \Big) \Big\} \\[2ex]
    	&\le \psi_1\Big(\frac{n  - p -1}{2}\Big) = \frac{2}{n-p-1}+O
    	\Big(\frac{1}{(n+p-1)^2}\Big) = o(1).
    \end{align*}
(note that  $p/n\to c<1$), which proves   \ref{A4}.  The conditions  \ref{A1} und \ref{A2} are obviously satisfied and the remaining inequality
\ref{A3} for the moments is a a consequence of Lemma A.7 and Theorem A.8. in \cite{tomdet2019}. Consequently, we obtain from
Theorem \ref{clt} the weak convergence
$$
S_n - \mathbb{E}[S_n] \to \mathcal{N}(0, \sigma^2)~,
$$
and  remains to calculate the representation of the expectation. For this purpose note that
    \begin{align}
        \mathbb{E}[S_n] &= \sum \limits_{i = 2}^{q_n} \sum \limits_{j = 1}^{p_i}
        \psi_{0}((n-p_i^*-j)/2) - \psi_{0}((n  - j)/2)    = B_{n}^{(1)} -  B_{n}^{2} \label{zweisummen}
       \end{align}
       where
         \begin{align}\label{bs}
        B_{n}^{(1)}   &= \sum  \limits_{i = n  - p}^{n - p_1-1} \psi_{0}(i/2)   ~,~
      B_{n}^{(2)}   =  \sum \limits_{i = 2}^{q_n} \sum \limits_{j = 1}^{p_i} \psi_{0}((n  - j)/2).
    \end{align}
Observing the expansion
$\psi_0(z) = \log(z) - \frac{1}{2z} + O(z^{-2})$ and \eqref{harmon}  we obtain
    \begin{align*}
           B_{n}^{(1)}     & = \sum  \limits_{i = n  - p}^{n - p_1-1}\Big\{ \log(i/2) - \frac{1}{i} + O(i^{-2}) \Big\}  \\[1ex]
        ={}& \log \Big(\frac{(n - p_1-1)!}{(n - p-1)!}\Big) - \log(2) (p - p_1)
        - \log\Big(\frac{n - p_1-1}{n - p-1}\Big) - \frac{p_1-p}{(n-p_1-1)(n-p-1)}
        \\[1ex] &{} +O((n-p_1-1)^{-2}) + O((n-p-1)^{-2}) \\[1ex]
        &= \log \Big(\frac{(n - p_1-1)!}{(n - p-1)!}\Big) - \log(2) (p - p_1)
        - \log\Big(\frac{n - p_1-1}{n - p-1}\Big) +o(1).
    \end{align*}
and
    \begin{align*}
            B_{n}^{(2)}    &= \sum \limits_{i = 2}^{q_n} \sum \limits_{j = 1}^{p_i} \psi_0((n  - j)/2) = \sum \limits_{i = 2}^{q_n} \sum \limits_{j = 1}^{p_i}\Big\{ \log(n - j) - \log(2) - \frac{1}{n  - j} + O((n - j)^{-2})\Big\} \\
        &= \Big\{\log\Big(\prod \limits_{i = 2}^{q_n} \frac{(n-1)!}{(n - p_i-1)!}\Big) - \log(2)(p - p_1) + \sum \limits_{i = 2}^\infty \log(1 - \lambda_i)\Big\} + o(1),
    \end{align*}
 where we used similar arguments as in the derivation of \eqref{a1} and \eqref{sum2}.
Combining these results with \eqref{zweisummen} finally yields
    \begin{align*}
        \mathbb{E}[S_n] & = \log\Big(\prod \limits_{i = 1}^{q_n}
        (n - p_i-1)!\Big) - \log((n - p - 1)!
        (n-1)!^{q_n - 1}) - \frac{\sigma^2}{2} + o(1).
    \end{align*}
Finally an application of  Stirling's formula
    \begin{align*}
    	\log(n!) = n \log n -n + \frac{1}{2} \log (2\pi n) +
    	\frac{1}{12  n} + O(n^{-3} )
	\end{align*}
   yields for the first term
    \begin{align*}
        &\log\Big(\prod \limits_{i = 1}^{q_n} (n - p_i - 1)!
        \Big ) - \log((n - p - 1)! (n-1)!^{q_n - 1}) \\[1ex]
        ={}&\sum \limits_{i = 1}^{q_n} \Big\{(n - p_i-1)
         \log(n - p_i-1) - (n - p_i - 1 )
         + \frac{\log(2\pi(n - p_i-1))}{2}
         + \frac{1}{12(n - p_i-1)} \\[1ex] &\qquad+ O((n - p_i-1)^{-3})
         \Big\}
         - \Big\{(n - p - 1)\log(n-p - 1)
        + (q_n - 1)(n-1)\log(n-1) \\[1ex] &\quad - (n - p - 1)  -(q_n-1)
        (n-1)  + \frac{\log(2\pi(n  - p - 1 )) +
        (q_n - 1)\log(2\pi (n-1))}{2} \\[1ex]
        & \qquad+ \frac{1}{12(n - p - 1)}  + \frac{q_n - 1}{12(n-1)}
         + O((n - p- 1)^{-3}) + q_n
          O((n-1)^{-3})\Big\} \\
        ={}&\sum \limits_{i = 1}^{q_n} \Big\{(n - p_i - 1)
        \log(1 - p_i/(n - 1)) + \frac{\log(1 - p_i/(n - 1))}{2} +
         \frac{1}{12(n - p_i - 1)} - \frac{1}{12(n-1)}\Big\} \\[1ex]
        & - \Big\{(n - p - 1)\log(1 - p/(n - 1)) ) +
        \frac{\log(1  - p/(n - 1))}{2} + \frac{1}{12(n - p -1)}  -
         \frac{1}{12(n-1)}\Big\}+ o(1)\\[1ex]
        ={}&\Big\{ \sum
        \limits_{i = 1}^{q_n} (n - p_i - 1) \log(1 - p_i/(n - 1))
        \Big\}   - (n - p - 1)\log(1 - p/(n - 1) ) + \frac{\sigma^2}{4} + o(1), \\[1ex]
    \end{align*}
   and the representation for the expectation in \eqref{b2} follows, completes the proof Theorem \ref{thm1}.

\subsection{Proof of Theorem \ref{thm2} }

    Define  $
    	\tilde{U}_n = \Lambda_n^{2/n}
$
    and note that under the null hypothesis, the distribution of $\tilde{U}_n$ is given by a product of independent Beta distributions [see \cite{anderson2003}], that is,
    \begin{align*}
    	\tilde{U}_n = \prod\limits_{i=1}^p U_{i},
    \end{align*}
    where the random variables
    \begin{align*}
    U_{i} \sim \beta \lb \frac{1}{2} \lb n - q_n + 1 - i \rb, \frac{1}{2} q_{1,n} \rb .
    \end{align*}
    Now consider the transformation
    	\begin{align*}
    	S_n := \frac{2}{n} \log \Lambda_n = \sum\limits_{i=1}^p \log (U_{i}),
    	\end{align*}
   then the assertion follows from
    \begin{align*}
     S_n - s_n \xrightarrow{{\cal D}_{H_0}} \mathcal{N}(0, \sigma^2),
    \end{align*}
    where $s_n$ and $\sigma^2$ are defined in \eqref{c1} and \eqref{c2}.  In order to prove the asymptotic normality of $S_n - s_n$, we show that the conditions \ref{A1}-\ref{A5} hold, beginning with  a derivation of the variance.
    \begin{align*}
		\sum\limits_{i=1}^{p} \Var ( \log (U_i))
		=& \sum\limits_{i=1}^p \Big\{ \psi_1 \Big( \frac{1}{2}  ( n - q_n + 1 - i ) \Big) - \psi_1 \Big( \frac{1}{2}  (  n - q_{2,n} + 1 - i ) \Big) \Big\} \\[1ex]
		=& \sum\limits_{i=1}^p \Big\{ \frac{2}{n - q_n + 1 - i} + \mathcal{O} \Big( \frac{1}{( n - q_n + 1 - i )^2} \Big)
		- \frac{2}{n - q_{2,n} + 1 - i} \\[1ex] &+ \mathcal{O} \Big( \frac{1}{( n - q_{2,n} + 1 - i )^2} \Big) \Big\} \\[1ex]
		&= 2 \Big\{ \sum\limits_{i=n - q_n + 1 - p}^{n - q_n} \frac{1}{i} - \sum\limits_{i= n - q_{2,n} + 1 - p}^{n - q_{2,n}} \frac{1}{i} \Big\} + o(1) \\[1ex]
		&= 2 \Big\{ \log \ \Big( \frac{ n - q_n}{n - q_n - p} \Big) - \log \Big( \frac{n - q_{2,n}}{n - q_{2,n}- p} \Big) \Big\} + \mathcal{O} \Big( \frac{1}{n - q_n - p} \Big) \\[1ex]
		&= 2 \Big\{ \log \Big( \frac{1}{1 - y_2} \Big) - \log \Big(
\frac{y_1 + y_2}{y_1 + y_2 - y_1 y_2} \Big) \Big\} + o(1).
	\end{align*}
   Regarding the error terms, note that  the assumptions of Theorem \ref{thm2} imply $n - q_n - p \to \infty$ and $\frac{p}{n-q_n-p} = \mathcal{O}(1)$. We continue  expanding the expected value
\begin{align*}
		\sum\limits_{i=1}^p \E ( \log (U_i )) =& \sum\limits_{i=1}^p \Big\{ \psi_0\Big( \frac{1}{2} \lb n - q_n + 1 - i  \rb) \Big) - \psi_0 \Big( \frac{1}{2} \lb n - q_{2,n} + 1 - i \rb \Big) \Big\} \\[1ex]
		=& \sum\limits_{i=1}^p \Big\{ \log\Big( \frac{1}{2} \lb n - q_n + 1 - i\rb \Big) - \frac{1}{n - q_n + 1 - i} \\[1ex]
		& - \Big[ \log \Big( \frac{1}{2} \lb n - q_{2,n} + 1 - i \rb \Big) - \frac{1}{n - q_{2,n} + 1 - i} \Big] \Big\}+ o(1) \\[1ex]
		=& \sum\limits_{i= n - q_n + 1 - p}^{n - q_n} \Big\{ \log \Big( \frac{i}{2} \Big) - \frac{1}{i} \Big\}
		- \sum\limits_{i=n - q_{2,n} + 1 - p }^{n-q_{2,n}} \Big\{ \log \Big( \frac{i}{2} \Big) - \frac{1}{i} \Big\} + o(1) \\[1ex]
		=& \log \Big( \frac{ ( n - q_n)! }{(n - q_n - p)! }\Big)  - \log \Big( \frac{n - q_n}{n - q_n - p} \Big)
		- \log \Big( \frac{( n - q_{2,n})!}{(n- q_{2,n} - p )!} \Big) \\[1ex]
		&+ \log \Big( \frac{n-q_{2,n}}{n-q_{2,n}-p} \Big) + o(1) \\[1ex]
		=& \log \Big( \frac{( n - q_n - 1)!}{(n - q_n - p -1)!} \Big) - \log \Big( \frac{(n - q_{2,n} - 1)!}{(n - q_{2,n} - p - 1)! } \Big) + o(1) \\[1ex]
		=& ( n - q_n - 1) \log ( n - q_n - 1) + (n - q_{2,n} - p - 1) \log ( n - q_{2,n} - p - 1) \\[1ex]
		&- (n - q_n - p -1) \log ( n- q_n - p - 1) - (n-q_{2,n}-1) \log (n - q_{2,n} - 1) + \frac{\sigma^2}{4} + o(1).
	\end{align*}
	Furthermore, we have
	\begin{align*}
		\sup\limits_{1\leq i \leq p} \Var ( \log (V_i)) &= \sup\limits_{1\leq i \leq p} \Big\{ \psi_1 \Big( \frac{1}{2} \lb n - q_n + 1 - i \rb \Big)
		- \psi_1\Big( \frac{1}{2} \lb n - q_{2,n} + 1 - i\rb \Big) \Big\} \\[1ex]
		& \leq \psi_1 \Big( \frac{1}{2} \lb n - q_n + 1 - p \rb \Big) = \mathcal{O} \Big( \frac{2}{n - q_n + 1 - p} \Big) =o(1),
	\end{align*}
	which is condition \ref{A4}. Obviously, \ref{A1} and \ref{A2} are also satisfied. The inequality for the moments in \ref{A3} follows from Lemma A.7 and Theorem A.8
	in \cite{tomdet2019}. Therefore, all conditions \ref{A1}-\ref{A5} are satisfied and the assertion follows from Theorem \ref{clt}.

\bigskip\bigskip

{\bf Acknowledgements.} 
The authors would like to thank
M. Stein who typed this manuscript with considerable technical expertise.
The work of H. Dette and  N.  D\"ornemann was partially supported by the Deutsche
Forschungsgemeinschaft (DFG Research Unit 1735, DE 502/26-2, RTG 2131). 

\setlength{\bibsep}{1pt}
\begin{small}

\end{small}

\end{document}